\documentclass[draft]{amsart}

\usepackage[letterpaper, left=2.5cm, right=2.5cm, top = 1in, bottom = 1in]{geometry}

\newtheorem{thm}{Theorem}[section]
\newtheorem{lem}{Lemma}[section]
\newtheorem*{thm*}{Theorem}

\newtheorem{corollary}{Corollary}[section]

\theoremstyle{definition}
\newtheorem{definition}{Definition}[section]
\newtheorem{example}{Example}[section]

\newcommand{\eA}{\mathcal{A}}
\newcommand{\iA}{\mathcal{A}^\ast}
\newcommand{\Or}{\mathcal{O}}
\newcommand{\M}{\mathcal{M}}
\DeclareMathOperator{\indeg}{indeg}
\DeclareMathOperator{\Div}{Div}
\newcommand{\T}{\mathcal{T}}
\newcommand{\E}{\mathcal{E}}

\usepackage{mathtools}
\usepackage{subcaption}
\usepackage{multicol}
\usepackage{tikz}
\usepackage{hyperref}
\usepackage{graphicx}
\usetikzlibrary{decorations.markings, arrows}
\tikzset{
	mid arrow/.style={postaction={decorate, decoration={
				markings,
				mark=at position .6 with {\arrow{latex}},
	}}},
}
\tikzset{
	bi arrow/.style={postaction={decorate, 
			decoration={
				markings,
				mark=at position .3 with {\arrowreversed{latex}},
				mark=at position .7 with {\arrow{latex}}
	}}},
}

\tikzset{mark size = 10}

\title{Generalized break divisors and triangulations of Lawrence polytopes}
\author{Natasha Crepeau}
\address{Department of Mathematics, University of Washington, United States}
\email{ncrepeau@uw.edu}
\begin{document}

		\begin{abstract}
		Let $G$ be a connected graph of genus $g$. The Picard group of degree $g$, $\text{Pic}^g(G)$, is the set of equivalence classes of divisors on $G$ of degree $g$, where two divisors are equivalent if one can be reached from the other through a sequence of chip-firing moves. We construct sets of representatives of the equivalence classes in  $\text{Pic}^g(G)$ by defining a function $I_G$ on the spanning trees of $G$ from a triangulation of the Lawrence polytope of the cographic matroid $\M^\ast(G)$. Additionally, such sets of representatives correspond to stability conditions on the nodal curve dual to the graph $G$. We show that $I_G$ that are constructed from regular triangulations of Lawrence polytope correspond to classical stability conditions, which are induced by generic real-valued divisors on $G$.
	\end{abstract}
	
	\maketitle
\section{Introduction}
\renewcommand*{\thethm}{\Alph{thm}}
Let $G = (V,E)$ be a loopless connected graph with vertices $V$ and edges $E$. Let $\mathcal{ST}(G)$ denote the set of all spanning trees of $G$. For a spanning tree $T \in \mathcal{ST}(G)$ and a choice of orientation $\Or_{\mathcal{E}}$ of the edges $\mathcal{E}$ not in $T$, the \textit{integral break divisor} $D(T,\Or_\mathcal{E})$ induced by $T$ and $\Or_{\mathcal{E}}$ is the degree $g = |E(G)| - |V(G)| + 1$ divisor on $G$ resulting from placing a chip at the head of each oriented edge in $\Or_{\mathcal{E}}$. Every degree $g$ divisor on $G$ is linearly equivalent to a unique integral break divisor and that the number of integral break divisors is the number of spanning trees of $G$ {\cite[Theorem 1.3]{integral}}. Therefore, the integral break divisors are a canonical choice of representatives of divisors of degree $g$.

In this paper, we study non-canonical choices of representatives of divisors of degree $g$. Consider some function $I_G$ that assigns a divisor $I_G(T)$ of degree $\deg(I_G(T)) = \sum_{v \in V(G)} I_G(T)(v) = 0$ to each spanning tree $T$ of $G$. Then, \text{a generalized break divisor} is a divisor of the form $\mathcal{D}(T, \Or_{\mathcal{E}}) + I_G(T)$. For a fixed function $I_G$, we denote the set of all generalized break divisors of the form $\mathcal{D}(T, \Or_{\mathcal{E}}) + I_G(T)$ by $\mathcal{BD}_{I_G}(G)$. For every function $I_G$, $|\mathcal{BD}_{I_G}(G)| \geq |\mathcal{ST}(G)|$, with $I_G \cong 0$ recovering the integral break divisors \cite[Theorem 3.5.1]{yuen2018}. We wish to describe all functions $I_G$ such that $|\mathcal{BD}_{I_G}(G)| = |\mathcal{ST}(G)|$; equivalently, we want to describe all $I_G$ such that every divisor of degree $g$ on $G$ is linearly equivalent to a unique divisor in $\mathcal{BD}_{I_G}(G)$.

To describe all such $I_G$, we turn to a family of bijections constructed by Ding \cite{ding2023framework} between the bases of a regular matroid $\M$, which in the graphic case are the spanning trees of the graph $G$, viewed as the bases of the graphic matroid $\M(G)$, and orientations of $\M$ up to \textit{circuit} and \textit{cocircuit reversal}, which in the graphic case are orientations of $G$ up to \textit{cycle} and \textit{cocycle reversal}. This family of bijections includes the Backman-Baker-Yuen (BBY) or ``geometric bijections" constructed in \cite{backman_baker_yuen_2019} and the Bernardi bijection restricted to spanning trees constructed in \cite{bernardi-bij}. These bijections arise from \textit{triangulations} and \textit{dissections} of the \textit{Lawrence polytopes} $\mathcal{P}$ and $\mathcal{P}^\ast$ corresponding respectively to the regular matroid $\M$ and its dual $\M^\ast$. When $\M$ is a graphic matroid, each triangulation and dissection of $\mathcal{P}$, respectively $\mathcal{P}^\ast$, corresponds to a choice of an external, respectively internal, orientation of each of the spanning trees of $G$, called \textit{triangulating} and \textit{dissecting atlases} \cite[Theorem 1.28]{ding2023framework}. Our first result uses triangulations of $\mathcal{P}^\ast$ to construct functions $I_G$ such that $\mathcal{BD}_{I_G}(G)$ is a set of representatives of divisors of degree $g$. Given an orientation $\Or$ of a graph $G$, the divisor $D_\Or$ is given by $D_\Or(v) = \indeg(v) - 1$ for all $v \in V(G)$. To construct the function $I_G$ from a triangulation of $\mathcal{P^\ast}$, we consider the
divisor $D_{-(\vec{T^\ast}^c)}$ induced by the internal orientation $\vec{T^\ast}$ in the corresponding triangulating internal atlas $\iA$. 

\begin{thm}
		For a graph $G$, consider the graphic matroid $\M(G)$. Recall that the bases of $\M(G)$ are the spanning trees $T \in \mathcal{ST}(G)$. Let $q$ be a vertex of $G$, and $\iA$ be a triangulating internal atlas. Then
		$\mathcal{BD}_{I_G^{\sigma^\ast}}(G)$ is a set of representatives of divisors of degree $g$, where 
		\begin{align*}
			I_G^{\sigma^\ast}: \mathcal{ST}(G) &\to \Div^0(G) \\
			T &\mapsto D_{-(\vec{T^\ast})^c} + (q). 
		\end{align*}
\end{thm}

Furthermore, atlases coming from triangulations of $\mathcal{P}$ (resp. $\mathcal{P}^\ast$) are in bijection with \textit{triangulating cycle} (resp. \textit{triangulating cocycle signatures}) of $G$ \cite[Theorem 1.16]{ding2023framework}. When the bijection is restricted to \textit{regular triangulations} of the Lawrence polytope, the corresponding signatures are the acyclic signatures used to define the BBY bijection \cite[Theorem 1.29]{ding2023framework}.  Regular triangulations of polytopes are connected by \textit{flips}, an operation that transforms one triangulation of a polytope into another by ``flipping" along a minimally dependent set of the vertices.

Another perspective on generalized break divisors is viewing them as \textit{stability conditions} on nodal curves, which have a dual graph as in, for example, \cite{pagani2022geometry}, \cite{pagani2023stability}, or \cite{viviani2023new}. A full definition is given in Section \ref{sec:stability conditions}. These stability conditions are in one-to-one correspondence with sets of generalized break divisors that are a set of representatives of divisors of degree $g$ {\cite[Theorem 1.20]{viviani2023new}}. In our last main result, we use our construction of $I_G^{\sigma^\ast}$, where $\sigma^\ast$ is the acyclic cocycle signature corresponding to some regular triangulation of $\mathcal{P^\ast}$, to construct a particular stability condition, called a \textit{classical} stability condition. 

\begin{thm}
	Assume $\sigma^\ast$ is an acyclic cocircuit signature of the graphic matroid $\M(G)$. Then, $\mathcal{BD}_{I_G^{\sigma^\ast}}$ is a classical stability condition, where $I_G^{\sigma^\ast} : \mathcal{ST}(G) \to {\Div}^0(G)$ is induced by $\sigma^\ast$. 
\end{thm}
\subsection*{Acknowledgements} The author was partially supported by NSF CAREER DMS-2044564. The author thanks Farbod Shokrieh for his guidance and advice during this project. Additionally, the author would like to thank Changxin Ding, Gaku Liu, Nicola Pagani, Lilla T\'othm\'er\'esz, Rhys Wells, and Cameron Wright for the helpful discussions.
\renewcommand*{\thethm}{\arabic{section}.\arabic{thm}}
\section{Preliminaries}
\subsection{Graphs and Matroids}
Let $G = (V,E)$ be a loopless graph. The \textit{genus} of $G$ is $g = |E(G)| - |V(G)| + 1$. A \textit{subgraph} of $G$ is a graph $G' = (V', E')$, where $V' \subseteq V(G)$ and $E' \subseteq E(G)$, and a subgraph $G'$ is \textit{spanning} if $V' = V(G)$. If $W \subseteq V(G)$, the \textit{induced subgraph} $G[W]$ is the subgraph of $G$ containing all edges between any vertices in $W$. If the induced subgraphs $G[W]$ and $G[W^c]$ are both connected, then we say $W$ is \textit{biconnected}. A \textit{simple cycle} is a walk along the graph $G$ that starts and ends at the same vertex without using any edge or any vertex twice, and a \textit{cut} is a subset of edges $E'$ such that $G \backslash E'$ is disconnected. A \textit{cocycle}, or \textit{bond}, is a minimal cut of a graph under inclusion. A \textit{spanning tree} $T$ of $G$ is a connected spanning subgraph of $G$ that does not have any cycles. The collection of all spanning trees of a graph $G$ is denoted by $\mathcal{ST}(G)$. Given a spanning tree $T$ and an edge $e \in G \backslash T$, the \textit{fundamental cycle} $C(T,e)$ is the unique cycle in $T \cup \{e\}$. If $e \in T$, $T \backslash e$ has two connected components with vertices $W$ and $W^c$. The \textit{fundamental cocycle} $C^\ast(T,e)$ is the unique cocycle partitioning $V(G)$ into $W$ and $W^c$. 

A \textit{divisor} $D$ on a graph $G$ is an assignment of integers $D(v)$ to the vertices of $G$. The \textit{degree} of a divisor $D$ is $\deg(D) = \sum_{v \in V(G)} D(v)$. We let $\Div^d(G)$ be the collection of all divisors of degree $d$ on $G$. A \textit{chip-firing move} takes a divisor $D$ to a divisor $D'$ by either \textit{firing at $v$}, where $v$ sends a chip to each vertex adjacent to it, or \textit{borrowing at $v$}, where $v$ receives a chip from each vertex adjacent to it.
If $D'$ is reached from $D$ through a sequence of chip-firing moves, we say that $D$ and $D'$ are \textit{linearly (or chip-firing) equivalent}. The equivalence class of a divisor $D$ is denoted by $[D]$. If $D$ is degree $k$, then $[D] \in \text{Pic}^k(G)$, the set of all divisors of degree $k$ of $G$ modulo chip-firing equivalence.

If $e = vw$ is an edge in $E(G)$, we can orient $e$ by assigning it a direction of towards $v$ or away from $v$. An \textit{orientation} $\Or$ of a graph is a choice of direction for every edge of $G$. The \textit{indegree} of $\Or$ at a vertex $v$ is the number of edges incident to $v$ pointing towards $v$. Given an orientation $\Or$, we can associate a divisor $D_\Or$ to it by setting $D_\Or = \sum_{v \in V(G)} (\text{indegree}(v) - 1)(v).$ The degree of $D_\Or$ is always $g-1$. If $D = D_\Or$ for some orientation $\Or$, we say that $D$ is \textit{orientable}. In addition, we will consider \textit{fourientations} of a graph, which were introduced in \cite{fourientations}. A fourientation generalizes orientations by allowing unoriented edges and edges oriented in both directions (bioriented), in addition to one-way oriented edges.
Given an orientation $\Or$ on a graph $G$, if $\Or$ has an oriented cycle $\vec{C}$ or cocycle $\vec{C^\ast}$, a \textit{cycle reversal} replaces $\vec{C}$ with $-\vec{C}$ in $\Or$, and a \textit{cocycle reversal} replaces $\vec{C^\ast}$ with $-\vec{C^\ast}$ in $\Or$. Gioan introduced the \textit{cycle-cocycle reversal system} in \cite{gioan}, where two orientations are considered to be equivalent if one can be reached from the other through a sequence of cycle and cocycle reversals. If $G$ is a connected graph, there are $|\mathcal{ST}(G)|$ such equivalence classes.
Cycle-cocycle reversal systems and chip-firing are connected, as shown in the following result. 
\begin{thm}[{\cite[Theorem 3.5]{BACKMAN2017655}}] \label{thm:reversal=chip-firing}
	Two orientations $\Or$ and $\Or'$ are equivalent in the cycle-cocycle reversal system if and only if $D_\Or$ is linearly equivalent to $D_{\Or'}$.
\end{thm}
For a graph $G$, consider its graphic matroid $\M(G)$, whose independent sets are the subgraphs of $G$ that do not contain a cycle. The \textit{circuits} of $\M(G)$ are the subsets corresponding to simple cycles of $G$. The matroid $\M(G)$ has a dual $\M^\ast(G)$, whose circuits correspond to the cocycles of the graph $G$. Graphic matroids are an example of \textit{regular matroids}, which are matroids represented by a totally unimodular matrix $M$. The size of matrix $M$ is $r \times n$, where $r$ is the rank of the matroid and $n = |E|$. For a graphic matroid $\M(G)$, $r(\M(G)) = |V(G)| - 1$, and $r(\M^\ast(G)) = g$, where $g$ is the genus of $G$. 

If $G$ is a directed graph, our matroid can have more structure. Let $G$ be a graph with an orientation $\Or$, and give every simple cycle $C$ of $G$ an orientation. An \textit{arc} is an oriented edge $\vec{e}$. Let $\vec{C}$ denote the choice of orientation for that cycle. Then, every cycle $C$ in the orientation $\Or$ has a positive part $C_+$, where the orientation of edge $e$ agrees with the arc $\vec{e} \in \vec{C}$, and a negative part $C_-$, where orientation of the edge $e$ is opposite of $\vec{e} \in \vec{C}$. We can write a signed cycle $\vec{C}$ as $\vec{C} = (C_+, C_-)$.

As is explained, for example, in \cite[Definition 3.2.1]{bjorner}, the collection of all signed cycles make $\M(G)$ an \textit{oriented matroid}. Picking an orientation for each of the cocircuits of $\M(G)$ would also yield an oriented matroid, and we denote signed cocircuits by $\vec{C^\ast}$. In \cite{ding2023framework}, the definition of fourientation is also extended to regular matroids, where a fourientation is a subset of the set of all arcs of $\M$; as in the graph case, $e \in E$ can be bioriented, unoriented, or one-way oriented. For a fourientation $\vec{F}$, $-\vec{F}$ is the fourientation with all one-way oriented arcs reversed, while $\vec{F}^c$ reverses one-way oriented arcs and switches unoriented and bioriented arcs. 

Oriented regular matroids also have a \textit{circuit-cocircuit reversal system}, introduced in \cite{gioan-circuit-cocircuit}, where an orientation $\Or$ of $\M$ is equivalent to another orientation $\Or'$ if $\Or'$ is reachable from $\Or$ through a sequence of circuit and cocircuit reversals. The number of circuit-cocircuit reversal classes is the number of bases of the matroid. The following lemma is helpful when studying these classes.
\begin{lem}[{\cite[Lemma 2.7]{ding2023framework}}] \label{lem:disjoint-circuits-cocircuits}
	Let $\Or_1$ and $\Or_2$ be two orientations in the same circuit-cocircuit reversal class of $\M$. Then $\Or_2$ can be obtained by reversing signed circuits and signed cocircuits in $\Or_1$ that do not share any edges.
\end{lem}
Choosing an orientation $\vec{C}$ for every circuit $C$ forms a \textit{circuit signature} $\sigma$, where $\sigma(C)=\vec{C}$. A \textit{cocircuit signature} $\sigma^\ast$ is defined similarly. We can identify $\sigma(C)$ with a vector in $\mathbb{R}^{|E|}$ with respect to a reference orientation of $\M$ by having the entry corresponding to $e$ be zero if $e \notin C$, $1$ if $\vec{e} \in \vec{C}$, and $-1$ if $-\vec{e} \in \vec{C}$. A circuit signature is \textit{acyclic} if for nonnegative reals $a_C$ such that $\sum_{C} a_C\sigma(C) = 0$, we must have $a_C = 0$ for all circuits $C$. Acyclic cocircuit signatures are defined similarly.
\subsection{Polytopes and Triangulations}
Let $\M$ be a loopless regular matroid that is represented by the totally unimodular $r \times n$ matrix $M$, written as $M = [I_r | L]$. Then, the dual matroid $\M^\ast$ can be represented by the matrix $M^\ast = [-L^T | I_{n-r}]$, as shown in Theorem 2.2.8 of \cite{oxley}. The \textit{Lawrence matrix} of $\M$ is given by
$$\Lambda(\M) = \begin{pmatrix}
	M_{r\times n} & \mathbf{0} \\
	I_{n\times n} & I_{n\times n}
\end{pmatrix}.$$
We label the columns of the Lawrence matrix as $P_1,...,P_n,P_{-1},...,P_{-n}$. The \textit{Lawrence polytope} $\mathcal{P} \subset \mathbb{R}^{n+r}$ of $\M$ is the convex hull of the points $P_1,...,P_n,P_{-1},...,P_{-n}$. If we replace $M$ with the matrix $M^\ast$, we get the Lawrence polytope $\mathcal{P}^\ast$ of the dual matroid $\M^\ast$. As shown in Proposition 4.4 of \cite{ding2023framework}, when $\M$ is loopless, the vertices of $\mathcal{P}$ are exactly the columns of $\Lambda(\M)$.

A \textit{triangulation} $\T$ of a polytope $\mathcal{P}$ is a collection of maximal simplices of $\mathcal{P}$ such that their union is $\mathcal{P}$ and any two distinct maximal simplices intersect at a (possibly empty) common face. If the second condition is replaced with any two distinct maximal simplices having a disjoint relative interior but not necessarily intersecting at a common face, then the collection of maximal simplices form a \textit{dissection} of $\mathcal{P}$. 

One special kind of triangulation is a \textit{regular triangulation}. We recall the following definition from \cite{santos2006geometric}. Let $V$ be the vertices of a polytope $\mathcal{P} \subset \mathbb{R}^d$, and let $\omega: V \to \mathbb{R}$ be a function that lifts $V$ to $\mathbb{R}^{d+1}$ by sending $V$ to
\[V_\omega \coloneqq \{(v, \omega(v)) : v \in V \}.\]
A \textit{lower facet} of $\text{conv}(V_\omega)$ is a facet whose supporting hyperplane lies below the interior of $\text{conv}(V_\omega)$. Then, project the lower facets of $\text{conv}(V_\omega)$ back onto $\mathcal{P}$. If all of the lower facets are simplices, then their projections form the triangulation $\mathcal{T}_\omega$ of $\mathcal{P}$, which we call a regular triangulation. 

Starting from one triangulation of $\mathcal{P}$, we can find another through a \textit{flip operation}. A circuit $C$ is a minimally affinely dependent set of the vertices $V$ of $\mathcal{P}$; in other words, if $C = \{c_1,...,c_k\}$, we have $\lambda_1c_1 + \cdots + \lambda_k c_k = 0$ with $\sum_i \lambda_i = 0$, and all $\lambda_i$ are nonzero, then $C$ is a circuit. Therefore, we can view $C$ as a signed circuit $C = (C_+, C_{-})$, where $C_+ = \{c_i : \lambda_i > 0\}$ and $C_{-} = \{c_i: \lambda_i < 0\}$. The circuit $C$ (as a point configuration) has two triangulations: $\mathcal{T}_+^C = \{C \backslash \{c_i\} : c_i \in C_+\}$ and $\mathcal{T}_{-}^C = \{C \backslash \{c_i\} : c_i \in C_{-}\}$. We also define the \textit{link} of $\tau \subseteq V$ with respect to a triangulation $\T$ as
\[ \text{link}_\mathcal{T}(\tau) = \{p \subseteq V: p \cap \tau = \emptyset, p \cup \tau \in \mathcal{T}\}.\]
\begin{definition}[{\cite[Definition 1.14]{santos2006geometric}}]
	Let $\mathcal{T}_1$ be a triangulation of the polytope $\mathcal{P}$. Assume $\mathcal{T}_1$ contains $\mathcal{T}_+^C$ for a circuit $C$ of the vertices $V$ of $\mathcal{P}$. Suppose all cells $\tau \in \mathcal{T}_+^C$ have the same link $L$ in $\mathcal{T}_1$. 
	The circuit $C$ \textit{supports a flip} in $\mathcal{T}_1$ and the triangulation $\mathcal{T}_2$ of $\mathcal{P}$ is obtained from $\mathcal{T}_1$ by this flip:
	$$\mathcal{T}_2 \coloneqq \mathcal{T}_1 \backslash \{ p \cup \tau : p \in L, \tau \in \mathcal{T}_+^C\} \cup \{p \cup \tau : p \in L, \tau \in \mathcal{T}_{-}^C\}$$
\end{definition}
The \textit{flip graph} of triangulations of a polytope $\mathcal{P}$ is constructed by having the vertices correspond to triangulations, where two triangulations are connected by an edge if they differ by a flip. As is well-known, the flip graph of regular triangulations is the $1$-skeleton of the secondary polytope of $\mathcal{P}$, so it is connected.
\subsection{Generalized Break Divisors}
Let $T$ be a spanning tree of $G$, and denote the complement of $T$ in $G$ by $\mathcal{E}_T$. Let $\Or_{\mathcal{E}_T}$ be a choice of orientation for the edges in $\E_T$. leaving the edges in $T$ unoriented. By dropping a chip at the head of each oriented edge in $\Or_{\E_T}$, we get the \textit{(integral) break divisor} $\mathcal{D}(T, \Or_{\E_T})$. Theorem 1.1 of \cite{integral} shows that the integral break divisors form a complete set of chip-firing representatives.
     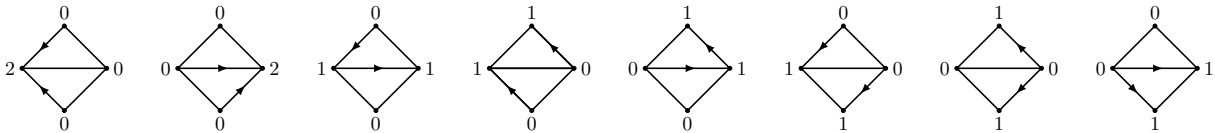
\begin{figure}[h!]
	\centering
	\scalebox{.75}{\begin{subfigure}[t]{.15\textwidth}
		\begin{tikzpicture}
			\draw[thick] (0,.75)--(.75,0)--(0,-.75);
			\draw[thick] (.75,0)--(-.75,0);
			\draw[thick, mid arrow] (0,.75) --(-.75,0);
			\draw[thick, mid arrow] (0,-.75) --(-.75,0);
			
			\filldraw[black] (0,.75) circle (1pt) node[anchor=south] {$0$};
			\filldraw[black] (.75,0) circle (1pt) node[anchor=west] {$0$};
			\filldraw[black] (0,-.75) circle (1pt) node[anchor=north] {$0$};
			\filldraw[black] (-.75,0) circle (1pt) node[anchor=east]{$2$};
		\end{tikzpicture}
	\end{subfigure}}
	\scalebox{.75}{
	\begin{subfigure}[t]{.15\textwidth}
		\begin{tikzpicture}
			\draw[thick] (.75,0)--(0,.75)--(-.75,0)--(0,-.75);
			\draw[thick, mid arrow] (-.75,0)--(.75,0);
			\draw[thick, mid arrow] (0,-.75)--(.75,0);
			\filldraw[black] (0,.75) circle (1pt) node[anchor=south] {$0$};
			\filldraw[black] (.75,0) circle (1pt) node[anchor=west] {$2$};
			\filldraw[black] (0,-.75) circle (1pt) node[anchor=north] {$0$};
			\filldraw[black] (-.75,0) circle (1pt) node[anchor=east]{$0$};
		\end{tikzpicture}     
	\end{subfigure}}
	\scalebox{.75}{
	\begin{subfigure}[t]{.15\textwidth}
		\begin{tikzpicture}
			\draw[thick] (-.75,0)--(0,-.75)--(.75,0)--(0,.75);
			\draw[thick, mid arrow] (0,.75) -- (-.75, 0);
			\draw[thick, mid arrow] (-.75,0)--(.75,0);
			\filldraw[black] (0,.75) circle (1pt) node[anchor=south] {$0$};
			\filldraw[black] (.75,0) circle (1pt) node[anchor=west] {$1$};
			\filldraw[black] (0,-.75) circle (1pt) node[anchor=north] {$0$};
			\filldraw[black] (-.75,0) circle (1pt)  node[anchor=east] {$1$};
		\end{tikzpicture}
	\end{subfigure}}
	\scalebox{.75}{
	\begin{subfigure}[t]{.15\textwidth}
		\begin{tikzpicture}
			\draw[thick] (0,-.75)--(.75,0)--(-.75,0);
			\draw[thick] (-.75,0)--(0,.75);
			\draw[thick, mid arrow] (0,-.75)--(-.75,0);
			\draw[thick, mid arrow] (.75,0)--(0,.75);
			\draw[thick] (0,.75)--(.75,0)--(-.75,0)--(0,-.75);
			
			\filldraw[black] (0,.75) circle (1pt) node[anchor=south] {$1$};
			\filldraw[black] (.75,0) circle (1pt) node[anchor=west] {$0$};
			\filldraw[black] (0,-.75) circle (1pt) node[anchor=north] {$0$};
			\filldraw[black] (-.75,0) circle (1pt) node[anchor=east] {$1$};
		\end{tikzpicture}
	\end{subfigure}}
	\scalebox{.75}{
	\begin{subfigure}[t]{.15\textwidth}
		\begin{tikzpicture}
			\draw[thick] (0,.75)--(-.75,0)--(0,-.75)--(.75,0);
			\draw[thick, mid arrow] (.75,0)--(0,.75);
			\draw[thick, mid arrow] (-.75,0)--(.75,0);
			\filldraw[black] (0,.75) circle (1pt) node[anchor=south] {$1$};
			\filldraw[black] (.75,0) circle (1pt) node[anchor=west] {$1$};
			\filldraw[black] (0,-.75) circle (1pt) node[anchor=north] {$0$};
			\filldraw[black] (-.75,0) circle (1pt) node[anchor=east]{$0$};
		\end{tikzpicture}
	\end{subfigure}}
	\scalebox{.75}{
	\begin{subfigure}{.15\textwidth}
		\begin{tikzpicture}
			\draw[thick] (0,.75)--(.75,0)--(-.75,0)--(0,-.75);
			\draw[thick, mid arrow] (0,.75)--(-.75,0);
			\draw[thick, mid arrow] (.75,0)--(0,-.75);
			\filldraw[black] (0,.75) circle (1pt) node[anchor=south] {$0$};
			\filldraw[black] (.75,0) circle (1pt) node[anchor=west] {$0$};
			\filldraw[black] (0,-.75) circle (1pt) node[anchor=north] {$1$};
			\filldraw[black] (-.75,0) circle (1pt) node[anchor=east] {$1$};
		\end{tikzpicture}
	\end{subfigure}}
	\scalebox{.75}{
	\begin{subfigure}[t]{.15\textwidth}
		\begin{tikzpicture}
			\draw[thick] (.75,0)--(-.75,0)--(0,-.75);
			\draw[thick] (-.75,0)--(0,.75);
			\draw[thick, mid arrow] (.75,0)--(0,.75);
			\draw[thick, mid arrow] (.75,0)--(0,-.75);
			\filldraw[black] (0,.75) circle (1pt) node[anchor = south] {$1$};
			\filldraw[black] (.75,0) circle (1pt) node[anchor=west] {$0$};
			\filldraw[black] (0,-.75) circle (1pt) node[anchor=north] {$1$};
			\filldraw[black] (-.75,0) circle (1pt) node[anchor=east]{$0$};
		\end{tikzpicture}
	\end{subfigure}}
	\scalebox{.75}{
	\begin{subfigure}[b]{.15\textwidth}
		\begin{tikzpicture}
			\draw[thick] (0,-.75)--(.75,0)--(0,.75)--(-.75,0);
			\draw[thick,mid arrow] (-.75,0)--(.75,0);
			\draw[thick, mid arrow] (-.75,0)--(0,-.75);
			\filldraw[black] (0,.75) circle (1pt) node[anchor=south] {$0$};
			\filldraw[black] (-.75,0) circle (1pt) node[anchor=east]{$0$};
			\filldraw[black] (0,-.75) circle (1pt) node[anchor=north] {$1$};
			\filldraw[black] (.75,0) circle (1pt)  node[anchor=west] {$1$};
		\end{tikzpicture}
	\end{subfigure}}
	
	\caption{Break divisors of the graph $G$, with an example $\Or_{\E_T}$ for each.}
	\label{fig:kite-break-divisors}
\end{figure}

In \cite{yuen2018}, Yuen, motivated by the work of Kass and Pagani in \cite{pagani-kass}, studied the notion of \textit{generalized break divisors}. Let $I_G: \mathcal{ST}(G) \to \Div^0(G)$, be any function that assigns a degree zero divisor to every spanning tree $T$ of $G$. 
\begin{definition}
	For a function $I_G: \mathcal{ST}(G) \to \Div^0(G)$, the set of \textit{generalized break divisors with respect to $I_G$} is defined as 
	\[\mathcal{BD}_{I_G}(G) \coloneqq \{\mathcal{D}(T, \Or_{\E_T}) + I_G(T) | \mathcal{D}(T, \Or_{\E_T}) \text{ is a break divisor of } G\}\]
\end{definition}
If $I_G \cong 0$, $\mathcal{BD}_{I_G}(G)$ is the set of integral break divisors of $G$, which we denote by $\mathcal{BD}_0(G)$. For a function $I_G$, the number of generalized break divisors is always bounded below by $|\mathcal{ST}(G)|$.
\begin{thm}{\cite[Theorem 3.5.1]{yuen2018}}
	For any function $I_G$, $|\mathcal{BD}_{I_G}(G)| \geq |\mathcal{ST}(G)|$. Equality holds when $I_G \cong 0$.
\end{thm}

\section{Stability conditions} \label{sec:stability conditions}
Given a nodal curve $X$ and its dual graph $G$, a question that one can ask is how to construct smoothable compactified Jacobians of $X$. In the fine case, the question was answered by Pagani and Tommasi in \cite{pagani2023stability} and by Viviani in \cite{viviani2023new} through the study of stability conditions on the graph $G$. For smoothable but not fine compacitified Jacobians, the question was answered by Fava, Pagani, and Viviani in \cite{fava-pagani-viviani2024}. 

\begin{definition}[{\cite[Definition 1.4]{viviani2023new}}]
	A \textit{nondegenerate $V$-stability condition} of degree $d \in \mathbb{Z}$ on $G$ is an assignment of integers 
	\[\mathfrak{n} = \{\mathfrak{n}_W : W \subset V(G) \text{ is biconnected and nontrivial }( W, W^c \neq \emptyset)\}\]
	satisfying the following properties: for any nontrivial biconnected $W \subset V(G)$, we have 
	\[\mathfrak{n}_W + \mathfrak{n}_{W^c} = d + 1 - |E(W, W^c)|;\]
	and any disjoint $W_1, W_2$ such that $W_1, W_2$ and $W_1 \cup W_2$ are all biconnected and nontrivial satisfy
	\[-1 \leq \mathfrak{n}_{W_1 \cup W_2} - \mathfrak{n}_{W_1} - \mathfrak{n}_{W_2} - |E(W_1, W_2)| \leq 0.\]
\end{definition}
Throughout this paper, all of our $V$-stability conditions will be nondegenerate and degree $g$, where $g$ is the genus of $G$. Degenerate $V$-stability conditions, corresponding to constructing smoothable but not fine compactified Jacobians, are defined and discussed in \cite{fava-pagani-viviani2024}.

There exists a smaller family of compactified Jacobians, called \textit{classical compactified Jacobians}, which have been constructed by several different authors such as Oda and Seshadri \cite{oda-seshadri}, Simpson \cite{simpson},  Caporaso \cite{caporaso}, and Esteves \cite{esteves}. In this paper, we work with the graph theoretic construction given in \cite{oda-seshadri}.

 A \textit{numerical polarization} $\phi: V(G) \to \mathbb{R}$ is a real-valued divisor on a graph $G$ such that the degree $\deg(\phi) = \sum_{v \in V(G)} \phi(v)$ is an integer $d$. A numerical polarization $\phi$ is \textit{generic} or \textit{nondegenerate} if $\sum_{v \in W} \phi(v) - \frac{|E(W, W^c)|}{2} = \phi_W - \frac{|E(W, W^c)|}{2} \notin \mathbb{Z}$ for all biconnected $W \subset V(G)$. Here $E(W, W^c)$ is the set of edges connecting a vertex in $W$ to a vertex in $W^c$.
Given a generic numerical polarization $\phi$, we can construct the $V$-stability condition $\mathfrak{n}(\phi)$, where each $\mathfrak{n}(\phi)_W = \left\lceil \phi_W - \frac{|E(W,W^c)|}{2} \right \rceil$. We say such a $V$-stability condition is \textit{classical}. 

As discussed in Section 1.3 of \cite{oda-seshadri} (see also \cite{pagani2023stability} and \cite{viviani2023new}), one can construct an infinite arrangement of hyperplanes in the affine space of numerical polarizations of degree $d$ supported on $V(G)$, $V^d(G)$, defined as
\[\mathcal{A}_G^d \coloneqq \left\{\phi_W - \frac{|E(W,W^c)|}{2} = n \right\}_{W \subset V(G), n \in \mathbb{Z}}\]
where $W$ again ranges over the biconnected subsets of $V(G)$. This yields a wall and chamber decomposition of $V^d(G)$ such that $\phi \in V^d(G)$ is only generic (or nondegenerate) if $\phi$ does not lie on any wall, and two generic polarizations $\phi$ and $\phi'$ are in the same chamber if and only if $\mathfrak{n}(\phi) = \mathfrak{n}(\phi')$. 

The following theorem connects nondegenerate $V$-stability conditions and sets of generalized break divisors.
\begin{thm}[{\cite[Theorem 1.20]{viviani2023new}}] \label{thm: bijection:vset-gbd}
	Let $G$ be a connected graph. There is a bijection
	\begin{align*}
		\{V\text{-stability conditions } \mathfrak{n}\} &\leftrightarrow \{\mathcal{BD}_{I_\mathfrak{n}}(G) \textit{ such that } |\mathcal{BD}_{I_\mathfrak{n}}(G)| = |\mathcal{ST}(G)| \}
	\end{align*}
	where $I_\mathfrak{n}$ is defined as follows: for $T \in \mathcal{ST}(G)$ and $e \in T$, $T \backslash e$ is the disjoint union of connected components $T[W_e]$ and $T[W_e^c]$, where $W_e$ and $W_e^c$ are biconnected subsets of $V(G)$. Then $I_\mathfrak{n}(T)$ is the unique element of $\Div^{0}(G)$ such that
	\[I_\mathfrak{n}(T)_{W_e} = \mathfrak{n}_{W_e} - g(G[W_e]) \text{ and } I_\mathfrak{n}(T)_{W_e^c} = \mathfrak{n}_{W_e^c} - g(G[W_e^c]),\]
	where $I_\mathfrak{n}(T)_{W_e} = -I_\mathfrak{n}(T)_{W_e^c} = \sum_{v \in W_e} I_\mathfrak{n}(T)(v)$. 
	Additionally, $|\mathcal{BD}_{I_G}(G)| = |\mathcal{ST}(G)|$ is equivalent to $\mathcal{BD}_{I_G}(G)$ being a complete set of chip-firing representatives.
\end{thm}

As a result, we refer to the $\mathcal{BD}_{I_G}(G)$ satisfying $|\mathcal{BD}_{I_G}(G)| = |\mathcal{ST}(G)|$ as stability conditions. If $\mathcal{BD}_{I_G}(G) = \mathcal{BD}_{I_{\mathfrak{n}(\phi)}}(G)$ for a generic numerical polarization $\phi$, we say $\mathcal{BD}_{I_G}(G)$ is classical.
\section{Bijections between signatures, atlases, and Lawrence polytopes}
Let $\M$ be a regular matroid. Let $B$ be a basis of $\mathcal{M}$. We call the edges in $B$ \textit{internal} and the edges not in $B$ \textit{external}. An \textit{externally (resp. internally) oriented basis} $\vec{B}$ (resp. $\vec{B^\ast}$) is a fourientation where all the internal (resp. external) edges are bioriented and all the external (resp. internal) edges are one-way oriented. An \textit{external (resp. internal) atlas} $\mathcal{A}$ (resp. $\iA$) of $\mathcal{M}$ is a collection of externally (resp. internally) oriented bases $\vec{B}$ such that each basis of $\mathcal{M}$ appears exactly once. External and internal atlases can be constructed from circuit and cocircuit signatures, respectively. Given a circuit signature $\sigma$, $\vec{B} \in \eA_\sigma$ is oriented such that if $e \notin B$, $e$ is oriented to be compatible with the signed circuit $\sigma(C(B, e))$. Similarly, given a cocircuit signature $\sigma^\ast$, $\vec{B^\ast} \in \iA_{\sigma^\ast}$ is oriented such that if $e \in B$, $e$ is oriented to be compatible with the signed cocircuit $\sigma^\ast(C^\ast(B,e))$.

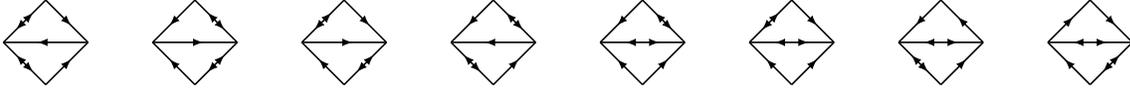
\begin{figure}[h!]
			\centering
			\scalebox{.75}{\begin{subfigure}[t]{.15\textwidth}
					\begin{tikzpicture}
							\draw[thick, mid arrow] (0,.75)--(.75,0);
							\draw[thick, mid arrow] (.75,0)--(-.75,0);
							\draw[thick, mid arrow] (0,-.75)--(.75,0);
							\draw[thick, bi arrow] (0, .75)--(-.75,0);
							\draw[thick, bi arrow] (-.75,0)--(0,-.75);
						\end{tikzpicture}
					\end{subfigure}}
			\scalebox{.75}{\begin{subfigure}[t]{.15\textwidth}
					\begin{tikzpicture}
							\draw[thick, mid arrow] (0,.75)--(-.75,0);
							\draw[thick, mid arrow] (0,-.75)--(-.75,0);
							\draw[thick, mid arrow] (-.75,0)--(.75,0);
							\draw[thick, bi arrow] (0, .75)--(.75,0);
							\draw[thick, bi arrow] (.75,0)--(0,-.75);
						\end{tikzpicture}
					\end{subfigure}}
				\scalebox{.75}{\begin{subfigure}[t]{.15\textwidth}
					\begin{tikzpicture}
							\draw[thick, mid arrow] (0,.75)--(.75,0);
							\draw[thick, mid arrow] (-.75,0)--(.75,0);
							\draw[thick, mid arrow] (0,-.75)--(-.75,0);
							\draw[thick, bi arrow] (0,.75)--(-.75,0);
							\draw[thick, bi arrow] (.75,0)--(0,-.75);
						\end{tikzpicture}
					\end{subfigure}}
			\scalebox{.75}{\begin{subfigure}[t]{.15\textwidth}
					\begin{tikzpicture}
							\draw[thick, mid arrow] (0,-.75)--(.75,0);
							\draw[thick, mid arrow] (0, .75)--(-.75,0);
							\draw[thick, mid arrow] (.75,0)--(-.75,0);
							\draw[thick, bi arrow] (0,-.75)--(-.75,0);
							\draw[thick, bi arrow] (.75,0)--(0,.75);
						\end{tikzpicture}
					\end{subfigure}}
			\scalebox{.75}{\begin{subfigure}[t]{.15\textwidth}
					\begin{tikzpicture}
							\draw[thick, mid arrow] (0,.75)--(-.75,0);
							\draw[thick, mid arrow] (0,-.75)--(.75,0);
							\draw[thick, mid arrow] (0,-.75)--(-.75,0);
							\draw[thick, bi arrow] (.75,0)--(0,.75);
							\draw[thick, bi arrow] (-.75,0)--(.75,0);
						\end{tikzpicture}
					\end{subfigure}}
			\scalebox{.75}{\begin{subfigure}[t]{.15\textwidth}
					\begin{tikzpicture}
							\draw[thick, mid arrow] (0,.75)--(.75,0);
							\draw[thick, mid arrow] (0,-.75)--(.75,0);
							\draw[thick, mid arrow] (0,-.75)--(-.75,0);
							\draw[thick, bi arrow] (-.75,0)--(0,.75);
							\draw[thick, bi arrow] (-.75,0)--(.75,0);
						\end{tikzpicture}
					\end{subfigure}}
			\scalebox{.75}{\begin{subfigure}[t]{.15\textwidth}
					\begin{tikzpicture}
							\draw[thick, mid arrow] (0,-.75)--(.75,0);
							\draw[thick, mid arrow] (.75,0)--(0,.75);
							\draw[thick, mid arrow] (0,.75)--(-.75,0);
							\draw[thick, bi arrow] (-.75,0)--(.75,0);
							\draw[thick, bi arrow] (-.75,0)--(0,-.75);
						\end{tikzpicture}
					\end{subfigure}}
			\scalebox{.75}{\begin{subfigure}[t]{.15\textwidth}
					\begin{tikzpicture}
							\draw[thick, mid arrow] (0,-.75)--(-.75,0);
							\draw[thick, mid arrow] (-.75,0) -- (0, .75);
							\draw[thick, mid arrow] (0,.75) -- (.75,0);
							\draw[thick, bi arrow] (-.75,0)--(.75,0);
							\draw[thick, bi arrow] (.75,0)--(0,-.75);
						\end{tikzpicture}         
					\end{subfigure}}
			\caption{Triangulating internal atlas of $\M(G)$}
			\label{fig:kite-internal atlas}
		\end{figure}

For the bijections in \cite{ding2023framework}, we want two types of atlases: dissecting and triangulating. Triangulating and dissecting atlases are constructed from triangulations and dissections of the corresponding Lawrence polytopes.
\begin{thm}[{\cite[Theorem 1.28]{ding2023framework}}]\label{thm: triang-to-atlas}
	Let $\mathcal{P}$ be the Lawrence polytope of the loopless matroid $\M$. The map
	\[\chi: \{\text{triangulations (resp. dissections) } \T \text{ of } \mathcal{P}\} \to \{\text{triangulating (resp. dissecting) external atlases of } \M \}\]
	is a bijection, where for a triangulation  (or dissection) $\T$ of $\mathcal{P}$, each simplex $\tau \in \T$ corresponds to an externally oriented base of $\M$. If we replace the Lawrence polytope $\mathcal{P}$ with $\mathcal{P}^\ast$, $\chi$ with $\chi^\ast$, and ``external" with ``internal" then the statement also holds.
\end{thm}
Circuit and cocircuit signatures can also be triangulating.
\begin{definition}[{\cite[Definition 1.14]{ding2023framework}}]
	A circuit signature $\sigma$ is said to be \textit{triangulating} if for any $\vec{B} \in \eA_\sigma$ and any signed circuit $\vec{C} \subseteq \vec{B}$, $\vec{C}$ is in the signature $\sigma$. A cocircuit signature $\sigma^\ast$ is said to be \textit{triangulating} if for any $\vec{B^\ast} \in \iA_{\sigma^\ast}$ and any signed cocircuit $\vec{C^\ast} \subseteq \vec{B^\ast}$, $\vec{C^\ast}$ is in the signature $\sigma^\ast$. 
\end{definition}
Lemma 3.4 in \cite{ding2023framework} states that acyclic signatures are always triangulating. 
\begin{thm}[{\cite[Theorem 1.16]{ding2023framework}}]\label{thm:triangsigs-triangatlas}
	The following maps are bijections.
	\begin{align*}
		\alpha : \{\text{triangulating circuit sig. of } \M\} &\to \{\text{triangulating external atlases of } \M\} \\
		\sigma &\mapsto \eA_\sigma\\
		\alpha^\ast : \{\text{triangulating cocircuit sig. of } \M\} &\to \{\text{triangulating internal atlases of } \M\} \\
		\sigma^\ast &\mapsto \iA_{\sigma^\ast}
	\end{align*}
\end{thm}
One example of a triangulating internal atlas is shown in Figure \ref{fig:kite-internal atlas}.
	
Restricting to regular triangulations gives another useful bijection.
\begin{thm}[{\cite[Theorem 1.29]{ding2023framework}}] \label{thm: acyclic-regular}
	The restriction of the bijection $\chi^{-1}\circ \alpha$ to the set of acyclic circuit signatures of $\M$ is bijective onto the set of regular triangulations of $\mathcal{P}$. Dually, the restriction of the bijection $(\chi^\ast)^{-1} \circ \alpha^\ast$ to the set of acyclic cocircuit signatures of $\M$ is bijective onto the set of regular triangulations of $\mathcal{P}^\ast$, where $\alpha$ and $\alpha^\ast$ are the bijections defined in Theorem \ref{thm:triangsigs-triangatlas}.
\end{thm}
The final bijection we present from \cite{ding2023framework} demonstrates the importance of dissecting and triangulating atlases.
\begin{thm}[{\cite[Theorem 1.7]{ding2023framework}}] \label{thm:base-to-cyclecocycle}
	Given a pair of dissecting atlases $(\eA, \iA)$ of a regular matroid $\mathcal{M}$, if at least one of the atlases is triangulating, then the map 
	\begin{align*}
		\bar{f}_{\mathcal{A}, \mathcal{A}^\ast}: \{\text{bases of } \mathcal{M} \}&\to \{\text{circuit-cocircuit reversal classes of } \mathcal{M} \}\\
		B & \mapsto [\vec{B} \cap \vec{B^\ast}]
	\end{align*}
	is bijective, where $[\vec{B} \cap \vec{B^\ast}]$ is the circuit-cocircuit reversal class containing the orientation $\vec{B} \cap \vec{B^\ast}$.
\end{thm}
\section{Triangulations and stability conditions}
Our main result is to connect the bijections in Theorems \ref{thm: triang-to-atlas} and \ref{thm:triangsigs-triangatlas} to stability conditions. Recall that given an orientation $\Or$, the divisor $D_\Or$ is given by $D_\Or(v) = \text{indegree}(v) - 1$ for all $v \in V(G)$. 

Recall that for a basis element $B$, the intersection $\vec{B} \cap \vec{B^\ast}$ of an external orientation and an internal orientation, which are fourientations of $B$, is some orientation of $B$. Additionally, recall that for the fourientation $\vec{B^\ast}$, $-(\vec{B^\ast})^c$ replaces the bioriented edges of $\vec{B^\ast}$ with unoriented edges. 
\begin{thm} \label{thm:triang-to-goodIG}
	For a graph $G$, consider the graphic matroid $\M(G)$. Recall that the bases of $\M(G)$ are the spanning trees $T \in \mathcal{ST}(G)$. Let $q$ be a vertex of $G$. Let $\eA$ and $\iA$ be a pair of dissecting atlases with $\iA = \iA_{\sigma^\ast}$ for some triangulating cocycle signature $\sigma^\ast$. 
	Then, the set $\{D_{\vec{T} \cap \vec{T^\ast}} + (q) | T \in \mathcal{ST}(G) \text{ and } \vec{T^\ast} \in \iA_{\sigma^\ast}\}$ is a set of pairwise linearly inequivalent generalized break divisors $\mathcal{BD}_{I_G^{\sigma^\ast}}(G)$, where 
	\begin{align*}
		I_G^{\sigma^\ast}: \mathcal{ST}(G) &\to \Div^0(G) \\
		T &\mapsto D_{-(\vec{T^\ast})^c} + (q) 
	\end{align*}
	and $\mathcal{BD}_{I_G^{\sigma^\ast}}(G)$ is a complete set of chip-firing representatives.
\end{thm}

\begin{proof}
	
	First, we define $I_G^{\sigma^\ast}(T) = D_{-(\vec{T^\ast})^c} + (q)$ for each spanning tree $T$ of $\mathcal{ST}(G)$. We have $\deg(I_G^{\sigma^\ast}(T)) = 0$, so $\mathcal{BD}_{I_G^{\sigma^\ast}}(G)$ is a set of generalized break divisors. 		The external atlas $\eA$ orients the external edges of each spanning tree $T$, which corresponds to placing a chip at the head of each one-way oriented edge, corresponding to some break divisor $\mathcal{D}(T, \Or_{\mathcal{E}_{T}})$, where $\mathcal{E}_T = G \backslash T$. Therefore, for each $D_{\vec{T} \cap \vec{T^\ast}} + (q)$, we have
	\begin{equation} \label{eqn: trees-indegree}
		D_{\vec{T} \cap \vec{T^\ast}} + (q) = \mathcal{D}(T,\Or_{\mathcal{E}_{T}}) + D_{-(\vec{T^\ast})^c} + (q) = D(T,\Or_{\mathcal{E}_{T}}) + I_G^{\sigma^\ast}(T)
	\end{equation}
	so the set $\{D_{\vec{T} \cap \vec{T^\ast}} + (q) | T \in \mathcal{ST}(G) \}$ is contained in $\mathcal{BD}_{I_G^{\sigma^\ast}}(G)$.
	
	Then, we wish to show that any two divisors $D_{\vec{T_i} \cap \vec{T_i^\ast}} + (q)$ and $D_{\vec{T_j} \cap \vec{T_j^\ast}} + (q)$ are linearly inequivalent. It suffices to show $D_{\vec{T_i} \cap \vec{T_i^\ast}}$ and $D_{\vec{T_j} \cap \vec{T_j^\ast}}$ are linearly inequivalent. The bijection $\bar{f}_{\eA, \iA_{\sigma^\ast}}$ from Theorem \ref{thm:base-to-cyclecocycle} sends $T_i$ and $T_j$ to $\vec{T_i} \cap \vec{T_i^\ast}$ and $\vec{T_j} \cap \vec{T_j^\ast}$, which are representatives of different cycle-cocycle reversal classes; therefore, the orientations $\vec{T_i} \cap \vec{T_i^\ast}$ and $\vec{T_j} \cap \vec{T_j^\ast}$ are not equivalent under cycle and cocycle reversal.  If $D_{\vec{T_i} \cap \vec{T_i^\ast}} \sim D_{\vec{T_j} \cap \vec{T_j^\ast}}$, then by Theorem \ref{thm:reversal=chip-firing} the orientations $\vec{T_i} \cap \vec{T_i^\ast}$ and $\vec{T_j} \cap \vec{T_j^\ast}$ are in the same cycle-cocycle reversal class, a contradiction.
	
	Therefore, $\{D_{\vec{T} \cap \vec{T^\ast}} + (q) | T \in \mathcal{ST}(G)\}$ are all pairwise linearly inequivalent. Because there are $|\mathcal{ST}(G)|$ divisors in $\{D_{\vec{T} \cap \vec{T^\ast}} + (q) | T \in \mathcal{ST}(G)\}$, the $D_{\vec{T} \cap \vec{T^\ast}} + (q)$ form a complete set of representatives of divisors in $\Div^g(G)$ under chip-firing.

	Finally, we show that $|\mathcal{BD}_{I_G^{\sigma^\ast}}(G)| = |\mathcal{ST}(G)|$.  Consider a divisor $D \in \mathcal{BD}_{I_G^{\sigma^\ast}}(G)$ which is not of the form $D_{\vec{T} \cap \vec{T^\ast}} + (q)$ for any spanning tree $T \in \mathcal{ST}(G)$. There exists a spanning tree $T_D \in \mathcal{ST}(G)$ such that $$D = I_G^{\sigma^\ast}(T_D) + \mathcal{D}(T_D, \Or_{\E_{T_D}})$$
	where the external orientation $\Or_{T_D}$ extending $\Or_{\E_{T_D}}$ to a fourientation is not the external orientation $\vec{T_D} \in \eA$. Since $\vec{T_D^\ast} \cap \Or_{T_D}$ is an orientation of $G$, there exists a spanning tree $T' \in \mathcal{ST}(G)$ such that $\vec{T_D^\ast} \cap \Or_{T_D}$ is a member of the cycle-cocycle reversal class $[\vec{T'} \cap \vec{T'^\ast}]$.
	
	First, note that $T' \neq T_D$: if $T' = T_D$, then by \eqref{eqn: trees-indegree}
	\[D = I_G^{\sigma^\ast}(T_D) + \mathcal{D}(T_D, \Or_{\E_{T_D}}) + (q) = I_G^{\sigma^\ast}(T') + \mathcal{D}(T', \Or_{\E_{T'}}) +(q) = D_{\vec{T'} \cap \vec{T'^\ast}} + (q).\]

	By assumption $\vec{T'} \cap \vec{T'^\ast}$ is not the same orientation as $\vec{T_D^\ast} \cap \Or_{T_D}$, so by Lemma \ref{lem:disjoint-circuits-cocircuits} the two orientations are related by reversals of disjoint signed cycles and cocycles. Cycle reversal does not change the associated divisor, so we focus on cocycles. Say $\vec{T_D^\ast} \cap \Or_{T_D}$ and $\vec{T'^\ast} \cap \vec{T'}$ differ by a sequence of disjoint cocycle reversals.  Let $\vec{C^\ast}$ be a signed cocycle of $\vec{T'} \cap \vec{T'^\ast}$ such that $-\vec{C^\ast}$ is in $\vec{T_D^\ast} \cap \Or_{T_D}$. Therefore, $\vec{C^\ast}$ is in the cocycle signature $\sigma^\ast$, since $\sigma^\ast$ is triangulating. However if $-\vec{C^\ast}$ is in $\vec{T_D^\ast} \cap \Or_{T_D}$, $-\vec{C^\ast}$ is in $\vec{T_D^\ast}$, which contradicts $\sigma^\ast$ being triangulating. 
	Therefore,  $\vec{T'} \cap \vec{T'^\ast}$ and $\vec{T_D^\ast} \cap \Or_{T_D}$ can only differ up to cycle reversal, and $D_{\vec{T_D^\ast} \cap \Or_{T_D}} = D_{\vec{T} \cap \vec{T^\ast}}$ for some spanning tree $T \in \mathcal{ST}(G)$.
\end{proof}
\begin{example} The first row of Figure \ref{fig:kite-gbd} details the construction of $I_G^{\sigma^\ast}$ from the triangulating internal atlas $\iA_{\sigma^\ast}$ in Figure \ref{fig:kite-internal atlas}. The divisors in the second row are exactly $\mathcal{BD}_{I_G^{\sigma^\ast}}(G)$, which are constructed by adding the break divisors $\mathcal{D}(T, \Or_{\E_T})$ in Figure \ref{fig:kite-break-divisors} to $I_G^{\sigma^\ast}(T)$. 
	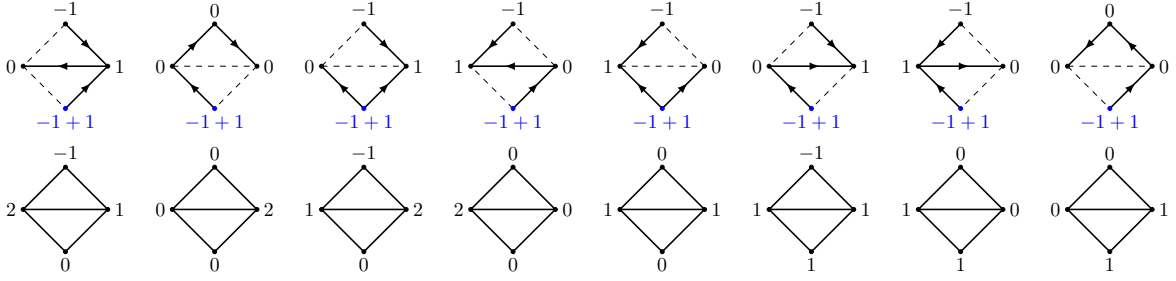
\begin{figure}[h!]
		\centering
		\scalebox{.75}{\begin{subfigure}[t]{.15\textwidth}
				\begin{tikzpicture}
					\draw[thick, mid arrow] (0,.75)--(.75,0);
					\draw[thick, mid arrow] (0,-.75)--(.75,0);
					\draw[thick, mid arrow] (.75,0)--(-.75,0);
					\draw[thin, dashed] (0, .75)--(-.75,0)--(0,-.75);
					
					\filldraw[black] (0,.75) circle (1pt) node[anchor=south] {$-1$};
					\filldraw[black] (-.75,0) circle (1pt) node[anchor=east] {$0$};
					\filldraw[blue] (0,-.75) circle (1pt) node[anchor=north] {$-1+1$};
					\filldraw[black] (.75,0) circle (1pt)  node[anchor=west] {$1$};
				\end{tikzpicture}
		\end{subfigure}}
		\scalebox{.75}{\begin{subfigure}[t]{.15\textwidth}
		\begin{tikzpicture}
			\draw[thick, mid arrow] (0,-.75)--(-.75,0);
			\draw[thick, mid arrow] (0,.75)--(.75,0);
			\draw[thick, mid arrow] (-.75,0)--(0,.75);
			\draw[thin, dashed] (-.75,0)--(.75,0);
			\draw[thin, dashed] (.75,0)--(0,-.75);
			\filldraw[black] (0,.75) circle (1pt) node[anchor=south] { $0$};
			\filldraw[black] (-.75,0) circle (1pt) node[anchor=east] {$0$};
			\filldraw[blue] (0,-.75) circle (1pt) node[anchor=north] {$-1+1$};
			\filldraw[black] (.75,0) circle (1pt)  node[anchor=west] {$0$};
		\end{tikzpicture}
\end{subfigure}}
		\scalebox{.75}{\begin{subfigure}[t]{.15\textwidth}
		\begin{tikzpicture}
			\draw[thick, mid arrow] (0,.75)--(.75,0);
			\draw[thin, dashed] (-.75,0)--(0,.75);
			\draw[thin, dashed] (-.75,0)--(.75,0);
			\draw[thick, mid arrow] (0, -.75) -- (.75,0);
			\draw[thick, mid arrow] (0, -.75) -- (-.75,0);
			\filldraw[black] (0,.75) circle (1pt) node[anchor=south] {$-1$};
			\filldraw[black] (-.75,0) circle (1pt) node[anchor=east] {$0$};
			\filldraw[blue] (0,-.75) circle (1pt) node[anchor=north] {$-1+1$};
			\filldraw[black] (.75,0) circle (1pt)  node[anchor=west] {$1$};
		\end{tikzpicture}
\end{subfigure}}
		\scalebox{.75}{\begin{subfigure}[t]{.15\textwidth}
				\begin{tikzpicture}
					\draw[thick, mid arrow] (0,-.75)--(.75,0);
					\draw[thick, mid arrow] (0,.75) -- (-.75,0);
					\draw[thick, mid arrow] (.75,0) -- (-.75,0);
					\draw[thin, dashed] (0,-.75)--(-.75,0);
					\draw[thin, dashed] (.75,0)--(0,.75);
					\filldraw[black] (0,.75) circle (1pt) node[anchor=south] {$-1$};
					\filldraw[black] (-.75,0) circle (1pt)  node[anchor=east] {$1$};
					\filldraw[blue] (0,-.75) circle (1pt) node[anchor=north] {$-1+1$};
					\filldraw[black] (.75,0) circle (1pt) node[anchor=west] {$0$};
				\end{tikzpicture}
		\end{subfigure}}
			\scalebox{.75}{\begin{subfigure}[t]{.15\textwidth}
			\begin{tikzpicture}
				\draw[thick, mid arrow] (0,.75)--(-.75,0);
				\draw[thick, mid arrow] (0, -.75) -- (.75,0);
				\draw[thick, mid arrow] (0, -.75) -- (-.75,0);
				\draw[thin, dashed] (.75,0)--(0,.75);
				\draw[thin, dashed] (-.75,0)--(.75,0);
				
				\filldraw[black] (0,.75) circle (1pt) node[anchor=south] {$-1$};
				\filldraw[black] (-.75,0) circle (1pt)  node[anchor=east] {$1$};
				\filldraw[blue] (0,-.75) circle (1pt) node[anchor=north] {$-1+1$};
				\filldraw[black] (.75,0) circle (1pt) node[anchor=west] {$0$};
			\end{tikzpicture}
	\end{subfigure}}
		\scalebox{.75}{\begin{subfigure}[t]{.15\textwidth}
		\begin{tikzpicture}
			\draw[thick, mid arrow] (0,.75)--(.75,0);
			\draw[thick, mid arrow] (0,-.75)--(-.75,0);
			\draw[thick, mid arrow] (-.75,0)--(.75,0);
			\draw[thin, dashed] (0,.75)--(-.75,0);
			\draw[thin, dashed] (.75,0)--(0,-.75);
			
			\filldraw[black] (0,.75) circle (1pt) node[anchor=south] {$-1$};
			\filldraw[black] (-.75,0) circle (1pt) node[anchor=east] {$0$};
			\filldraw[blue] (0,-.75) circle (1pt) node[anchor=north] {$-1+1$};
			\filldraw[black] (.75,0) circle (1pt)  node[anchor=west] {$1$};
		\end{tikzpicture}
\end{subfigure}}
		\scalebox{.75}{\begin{subfigure}[t]{.15\textwidth}
		\begin{tikzpicture}
			\draw[thick, mid arrow] (0,.75)--(-.75,0);
			\draw[thick, mid arrow] (-.75,0)--(.75,0);
			\draw[thick, mid arrow] (0,-.75)--(-.75,0);
			\draw[thin, dashed] (0, .75)--(.75,0)--(0,-.75);
			
			\filldraw[black] (0,.75) circle (1pt) node[anchor=south] {$-1$};
			\filldraw[black] (-.75,0) circle (1pt)  node[anchor=east] {$1$};
			\filldraw[blue] (0,-.75) circle (1pt) node[anchor=north] {$-1+1$};
			\filldraw[black] (.75,0) circle (1pt) node[anchor=west] {$0$};
		\end{tikzpicture}
\end{subfigure}}
		\scalebox{.75}{\begin{subfigure}[t]{.15\textwidth}
				\begin{tikzpicture}
					\draw[thick, mid arrow] (0,-.75)--(.75,0);
					\draw[thick, mid arrow] (0,.75)--(-.75,0);
					\draw[thick, mid arrow] (.75,0)--(0,.75);
					\draw[thin, dashed] (-.75,0)--(.75,0);
					\draw[thin, dashed] (-.75,0)--(0,-.75);
					\filldraw[black] (0,.75) circle (1pt) node[anchor=south] {$0$};
					\filldraw[black] (-.75,0) circle (1pt) node[anchor=east] {$0$};
					\filldraw[blue] (0,-.75) circle (1pt) node[anchor=north] {$-1+1$};
					\filldraw[black] (.75,0) circle (1pt)  node[anchor=west] {$0$};
				\end{tikzpicture}
		\end{subfigure}}

		\scalebox{.75}{\begin{subfigure}[t]{.15\textwidth}
				\begin{tikzpicture}
					\draw[thick] (0,.75)--(.75,0)--(0,-.75);
					\draw[thick] (.75,0)--(-.75,0);
					\draw[thick] (0, .75)--(-.75,0)--(0,-.75);
					
					\filldraw[black] (0,.75) circle (1pt) node[anchor=south] {$-1$};
					\filldraw[black] (-.75,0) circle (1pt) node[anchor=east] {$2$};
					\filldraw[black] (0,-.75) circle (1pt) node[anchor=north] {$0$};
					\filldraw[black] (.75,0) circle (1pt)  node[anchor=west] {$1$};
				\end{tikzpicture}
		\end{subfigure}}
		\scalebox{.75}{\begin{subfigure}[t]{.15\textwidth}
		\begin{tikzpicture}
			\draw[thick] (-.75,0)--(0,-.75)--(.75,0)--(0,.75);
			\draw[thick] (-.75,0)--(0,.75);
			\draw[thick] (-.75,0)--(.75,0);
			
			\filldraw[black] (0,.75) circle (1pt) node[anchor=south] {$0$};
			\filldraw[black] (-.75,0) circle (1pt) node[anchor=east] {$0$};
			\filldraw[black] (0,-.75) circle (1pt) node[anchor=north] {$0$};
			\filldraw[black] (.75,0) circle (1pt)  node[anchor=west] {$2$};
		\end{tikzpicture}
\end{subfigure}}
		\scalebox{.75}{\begin{subfigure}[t]{.15\textwidth}
				\begin{tikzpicture}
					\draw[thick] (.75,0)--(-.75,0)--(0,-.75);
					\draw[thick] (-.75,0)--(0,.75);
					\draw[thick] (0, .75)--(.75,0)--(0,-.75);
					
					\filldraw[black] (0,.75) circle (1pt) node[anchor=south] {$-1$};
					\filldraw[black] (-.75,0) circle (1pt)  node[anchor=east] {$1$};
					\filldraw[black] (0,-.75) circle (1pt) node[anchor=north] {$0$};
					\filldraw[black] (.75,0) circle (1pt) node[anchor=west] {$2$};
				\end{tikzpicture}
		\end{subfigure}}
		\scalebox{.75}{\begin{subfigure}[t]{.15\textwidth}
		\begin{tikzpicture}
			\draw[thick] (0,-.75)--(.75,0)--(-.75,0);
			\draw[thick] (-.75,0)--(0,.75);
			\draw[thick] (0,-.75)--(-.75,0);
			\draw[thick] (.75,0)--(0,.75);
			\filldraw[black] (0,.75) circle (1pt) node[anchor=south] {$0$};
			\filldraw[black] (-.75,0) circle (1pt)  node[anchor=east] {$2$};
			\filldraw[black] (0,-.75) circle (1pt) node[anchor=north] {$0$};
			\filldraw[black] (.75,0) circle (1pt) node[anchor=west] {$0$};
		\end{tikzpicture}
\end{subfigure}}
		\scalebox{.75}{\begin{subfigure}[t]{.15\textwidth}			
				\begin{tikzpicture}
					\draw[thick] (0,.75)--(-.75,0)--(0,-.75)--(.75,0);
					\draw[thick] (.75,0)--(0,.75);
					\draw[thick] (-.75,0)--(.75,0);
					
					\filldraw[black] (0,.75) circle (1pt) node[anchor=south] {$0$};
					\filldraw[black] (-.75,0) circle (1pt)  node[anchor=east] {$1$};
					\filldraw[black] (0,-.75) circle (1pt) node[anchor=north] {$0$};
					\filldraw[black] (.75,0) circle (1pt) node[anchor=west] {$1$};
				\end{tikzpicture}
		\end{subfigure}}
		\scalebox{.75}{\begin{subfigure}[t]{.15\textwidth}
		\begin{tikzpicture}
			\draw[thick] (0,.75)--(.75,0)--(-.75,0)--(0,-.75);
			\draw[thick] (0,.75)--(-.75,0);
			\draw[thick] (.75,0)--(0,-.75);
			
			\filldraw[black] (0,.75) circle (1pt) node[anchor=south] {$-1$};
			\filldraw[black] (-.75,0) circle (1pt) node[anchor=east] {$1$};
			\filldraw[black] (0,-.75) circle (1pt) node[anchor=north] {$1$};
			\filldraw[black] (.75,0) circle (1pt)  node[anchor=west] {$1$};
		\end{tikzpicture}
\end{subfigure}}
		\scalebox{.75}{\begin{subfigure}[t]{.15\textwidth}
		\begin{tikzpicture}
			\draw[thick] (.75,0)--(0,.75)--(-.75,0)--(0,-.75);
			\draw[thick] (-.75,0)--(.75,0);
			\draw[thick] (.75,0)--(0,-.75);
			\filldraw[black] (0,.75) circle (1pt) node[anchor=south] {$0$};
			\filldraw[black] (-.75,0) circle (1pt) node[anchor=east] {$1$};
			\filldraw[black] (0,-.75) circle (1pt) node[anchor=north] {$1$};
			\filldraw[black] (.75,0) circle (1pt)  node[anchor=west] {$0$};
		\end{tikzpicture}
\end{subfigure}}
		\scalebox{.75}{\begin{subfigure}[t]{.15\textwidth}
				\begin{tikzpicture}
					\draw[thick] (0,-.75)--(.75,0)--(0,.75)--(-.75,0);
					\draw[thick] (-.75,0)--(.75,0);
					\draw[thick] (-.75,0)--(0,-.75);
					\filldraw[black] (0,.75) circle (1pt) node[anchor=south] {$0$};
					\filldraw[black] (-.75,0) circle (1pt) node[anchor=east] {$0$};
					\filldraw[black] (0,-.75) circle (1pt) node[anchor=north] {$1$};
					\filldraw[black] (.75,0) circle (1pt)  node[anchor=west] {$1$};
				\end{tikzpicture}
		\end{subfigure}}
		\caption{Construction of $I_G^{\sigma^\ast}$ and the generalized break divisors in $\mathcal{BD}_{I_G^{\sigma^\ast}}$, with $q$ being the lower blue vertex.}
		\label{fig:kite-gbd}
	\end{figure}
\end{example}
Let $\sigma_1^\ast$ and $\sigma_2^\ast$ be two triangulating cocycle signatures which, by Theorem \ref{thm:triang-to-goodIG}, correspond to the functions $I_G^{\sigma_1^\ast}$ and $I_G^{\sigma_2^\ast}$. To better understand stability conditions of $G$, we wish to relate $I_G^{\sigma_1^\ast}$ and $I_G^{\sigma_2^\ast}$.
\begin{definition}
	Suppose $\sigma_1^\ast$ and $\sigma_2^\ast$ are two triangulating cocycle signatures. Then, we say they differ by a \textit{cocycle flip} if there exists a unique cocycle $C^\ast$ such that $\vec{C^\ast} \in \sigma_1^\ast$ and $-\vec{C^\ast} \in \sigma_2^\ast$.
\end{definition}
	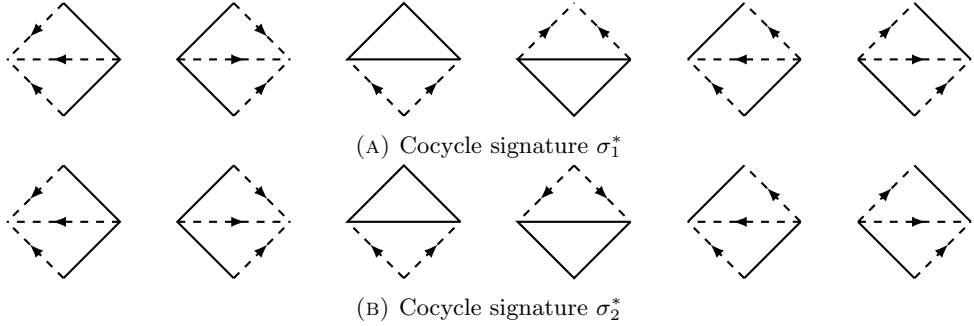
\begin{figure}[h!]
		\centering
		\begin{subfigure}{\textwidth}
			\centering
				\begin{tikzpicture}
					\draw[thick] (0,.75)--(.75,0)--(0, -.75);
					\draw[thick, dashed, mid arrow] (0, .75)--(-.75,0);
					\draw[thick, dashed, mid arrow] (.75,0)--(-.75,0);
					\draw[thick, dashed, mid arrow] (0,-.75)--(-.75,0);
				\end{tikzpicture}
			\hspace{.5cm}
				\begin{tikzpicture}
					\draw[thick] (0,.75)--(-.75,0)--(0, -.75);
					\draw[thick, dashed, mid arrow] (0, .75)--(.75,0);
					\draw[thick, dashed, mid arrow] (-.75,0)--(.75,0);
					\draw[thick, dashed, mid arrow] (0,-.75)--(.75,0);
				\end{tikzpicture}
			\hspace{.5cm}
				\begin{tikzpicture}
					\draw[thick] (.75,0)--(-.75,0)--(0,.75)--(.75,0);
					\draw[thick, dashed, mid arrow] (0,-.75)--(.75,0);
					\draw[thick, dashed, mid arrow] (0,-.75)--(-.75,0);   
				\end{tikzpicture}
			\hspace{.5cm}
				\begin{tikzpicture}
					\draw[thick] (.75,0)--(-.75,0)--(0,-.75)--(.75,0);
					\draw[thick, dashed, mid arrow] (.75,0)--(0,.75);
					\draw[thick, dashed, mid arrow] (-.75,0)--(0,.75);
				\end{tikzpicture}
			\hspace{.5cm}
				\begin{tikzpicture}
					\draw[thick] (0,.75)--(-.75,0);
					\draw[thick] (.75,0)--(0,-.75);
					\draw[thick, dashed, mid arrow] (.75,0)--(0,.75);
					\draw[thick, dashed, mid arrow] (0,-.75)--(-.75,0);
					\draw[thick, dashed, mid arrow] (.75,0)--(-.75,0);
				\end{tikzpicture}
			\hspace{.5cm}
				\begin{tikzpicture}
					\draw[thick] (0,-.75)--(-.75,0);
					\draw[thick] (.75,0)--(0,.75);
					\draw[thick, dashed, mid arrow] (-.75,0)--(0,.75);
					\draw[thick, dashed, mid arrow] (0,-.75)--(.75,0);
					\draw[thick, dashed, mid arrow] (-.75,0)--(.75,0);
				\end{tikzpicture}

			\subcaption{Cocycle signature $\sigma_1^\ast$}
		\end{subfigure}
		
		\begin{subfigure}{\textwidth}
			\centering
				\begin{tikzpicture}
					\draw[thick] (0,.75)--(.75,0)--(0, -.75);
					\draw[thick, dashed, mid arrow] (0, .75)--(-.75,0);
					\draw[thick, dashed, mid arrow] (.75,0)--(-.75,0);
					\draw[thick, dashed, mid arrow] (0,-.75)--(-.75,0);
				\end{tikzpicture}
			\hspace{.5cm}
				\begin{tikzpicture}
					\draw[thick] (0,.75)--(-.75,0)--(0, -.75);
					\draw[thick, dashed, mid arrow] (0, .75)--(.75,0);
					\draw[thick, dashed, mid arrow] (-.75,0)--(.75,0);
					\draw[thick, dashed, mid arrow] (0,-.75)--(.75,0);
				\end{tikzpicture}
			\hspace{.5cm}
				\begin{tikzpicture}
					\draw[thick] (.75,0)--(-.75,0)--(0,.75)--(.75,0);
					\draw[thick, dashed, mid arrow] (0,-.75)--(.75,0);
					\draw[thick, dashed, mid arrow] (0,-.75)--(-.75,0);   
				\end{tikzpicture}
			\hspace{.5cm}			
				\begin{tikzpicture}
					\draw[thick] (.75,0)--(-.75,0)--(0,-.75)--(.75,0);
					\draw[thick, dashed, mid arrow] (0,.75)--(.75,0);
					\draw[thick, dashed, mid arrow] (0,.75)--(-.75,0);   
				\end{tikzpicture}
			\hspace{.5cm}
				\begin{tikzpicture}
					\draw[thick] (0,.75)--(-.75,0);
					\draw[thick] (.75,0)--(0,-.75);
					\draw[thick, dashed, mid arrow] (.75,0)--(0,.75);
					\draw[thick, dashed, mid arrow] (0,-.75)--(-.75,0);
					\draw[thick, dashed, mid arrow] (.75,0)--(-.75,0);
				\end{tikzpicture}
			\hspace{.5cm}
				\begin{tikzpicture}
					\draw[thick] (0,-.75)--(-.75,0);
					\draw[thick] (.75,0)--(0,.75);
					\draw[thick, dashed, mid arrow] (-.75,0)--(0,.75);
					\draw[thick, dashed, mid arrow] (0,-.75)--(.75,0);
					\draw[thick, dashed, mid arrow] (-.75,0)--(.75,0);
				\end{tikzpicture}
			\subcaption{Cocycle signature $\sigma_2^\ast$}
		\end{subfigure}
		\caption{Two triangulating cocycle signatures that differ by a cocycle flip in the fourth cocycle. }
		\label{fig:kite-cocircuit-reachable}
	\end{figure}
The following lemma follows immediately from the definition of being triangulating. 
\begin{lem}
	If $\sigma_1^\ast$ and $\sigma_2^\ast$ are two triangulating cocycle signatures that differ by a cocycle flip at $C^\ast$, then, for each base $B$ such that $\vec{C^\ast} \subseteq \vec{B^\ast}$, $B$ only has one internal edge in $C^\ast$.
\end{lem}
In fact, using cocycle flips to move between triangulating cocycle signatures is the same as using flips to move between triangulations of the Lawrence polytope $\mathcal{P}^\ast$.
\begin{lem} \label{lem:cocycleflips_are_bistellar}
	Let $\mathcal{T}_1$ and $\mathcal{T}_2$ be triangulations of the Lawrence polytope $\mathcal{P^\ast}$. Let $\sigma_1^\ast$ be the cocycle signature such that $\chi^\ast(\mathcal{T}_1) = \iA_{\sigma_1^\ast}$, and define $\sigma_2^\ast$ similarly. Then, $\mathcal{T}_1$ and $\mathcal{T}_2$ differ by a flip if and only if $\sigma_1^\ast$ and $\sigma_2^\ast$ differ by a cocycle flip. 
\end{lem}
\begin{proof}
	First, assume that $\T_1$ and $\T_2$ differ by a flip. Therefore, there exists a circuit $C$ of the dual matroid $\M^\ast$ such that (without loss of generality) $\T_+^C$ is contained in $\T_1$ and $\text{link}_{\T_1}(\tau)$ is the same for every cell $\tau$ in $\T_+^C$. 
	
	We then note that $\T_2$ can be defined as
	$$\T_2 \coloneqq \mathcal{T}_1 \backslash \{ p \cup \tau : p \in \text{link}_{\T_1}(\tau), \tau \in \mathcal{T}_+^C\} \cup \{p \cup \tau : p \in \text{link}_{\T_1}(\tau), \tau \in \mathcal{T}_{-}^C\}. $$
	
	By duality, the signed circuit $C$ of $\M^\ast$ corresponds to the signed cocircuit $C^\ast = (C_+^\ast, C_{-}^\ast)$ of $\M(G)$. Recall the map $\chi^\ast$ from Theorem \ref{thm: triang-to-atlas}. For each $\tau \in \T_+^C$, the image of $\chi^\ast(\tau)$ is a fourientation of the subgraph formed by the edges of the cocycle $C^\ast$, where every edge is bioriented except for a single one-way oriented edge. Call that edge $e_\tau$. For each $\tau$, that edge $e_\tau$ must have  the same orientation, since each $\tau$ is constructed by $c_i \in C_+$ from $C$. We give the cocycle this orientation and denote it by $\vec{C^\ast}$.
	
	Since $\sigma^\ast_1$ is triangulating, we must have $\vec{C}^\ast = \sigma_1^\ast(C^\ast)$. The bistellar flip replaces $C_+$ with $C_-$, and therefore $\sigma_2^\ast(C^\ast) = -\vec{C^\ast}$, as desired.
	
	Now, assume that $\sigma_1^\ast$ and $\sigma_2^\ast$ differ by a cocycle flip. Let $C^\ast$ be that cocycle, and let $C$ be the circuit dual to it in $\M^\ast$. Because $\sigma_1^\ast$ is triangulating, we assume without loss of generality that $C_+$ is contained in $\T_1$. Let $L = \{ (\chi^\ast)^{-1}(\vec{B^\ast}) \backslash \tau | \vec{C^\ast} \subseteq \vec{B^\ast}, \vec{B^\ast} \in \iA_{\sigma_1^\ast} \}$. We wish to show that $L = \text{link}_{\T_1}(\tau)$ for every $\tau \in \T_{+}^{C}$.
	
	If $p \in \text{link}_{\T_1}(\tau)$, we have $p \cup \tau \in \T_1$ by definition of the link. Therefore, $\chi^\ast(p\cup \tau) = \vec{B^\ast}$ by Theorem \ref{thm: triang-to-atlas}. We then have $\vec{C^\ast} \subseteq \chi^\ast(\tau) \subseteq \vec{B^\ast}$, and therefore $p \in L$. For the reverse inclusion assume $\ell \in L$. Then, $\ell = (\chi^\ast)^{-1}(\vec{B^\ast}) \backslash \tau$ for some internally oriented base $\vec{B^\ast}$ of $\M(G)$. By definition $\ell \cap \tau = \emptyset$, and additionally $\chi^\ast(\ell) \cup \chi^\ast(\tau) = \vec{B^\ast}$ which is in the internal atlas $\iA_{\sigma_1^\ast}$. Therefore, $\ell \cup \tau \in \T_1$ and $\ell \in \text{link}_{\T_1}(\tau)$.
	
	Therefore, $C$ supports a bistellar flip in $\T_1$, and the result of that flip is $\T_2$.
\end{proof}
We can then consider the function $I^{\sigma_2^\ast}_G$ to be \textit{reachable} from $I_G^{\sigma_1^\ast}$ if $\sigma_2^\ast$ and $\sigma_1^\ast$ are connected by a sequence of cocycle flips. We recall that the flip graph restricted to regular triangulations is connected. 
\begin{corollary}\label{cor:reachable-regular}
	Let $\T_1$ and $\T_2$ be regular triangulations of $\mathcal{P}^\ast$, and $\sigma_1^\ast$ and $\sigma_2^\ast$ be the cocycle signatures such that $\chi^\ast(\T_1) = \iA_{\sigma_1^\ast}$ and $\chi^\ast(\T_2) = \iA_{\sigma_2^\ast}$. Then, $\sigma_1^\ast$ is reachable from $\sigma_2^\ast$. 
\end{corollary}
The first step towards constructing classical stability conditions from regular triangulations of $\mathcal{P}^\ast$ is to demonstrate that $\mathcal{BD}_0(G)$ is classical. We use the following lemma.
\begin{lem}\label{lem: subgraph_degrees}
	Let $W$ be a subset of $V(G)$. Then, $\sum_{v \in W} \deg(v) = |E(W, W^c)| + 2|E(G[W])|$.
\end{lem}
We then consider the generic polarization $\phi_{pcan}$ defined in Example 5.4 \cite{christ-payne-shen} or in Example 8.10 of \cite{pagani2023stability}, given by 
\[ \phi_{pcan}(v) = \frac{g(G) + |V(G)|}{2(g(G) + |V(G)|) - 2}\cdot \deg(v) - 1.\]
\begin{lem}\label{lem:n_w(pcan)}
	The $V$-stability condition $\mathfrak{n}$ associated to $\phi_{pcan}$ is
	\[\mathfrak{n}(\phi_{pcan}) = \{\mathfrak{n}(\phi_{pcan})_W = g(G[W]) \text{ for any biconnected and nontrivial }W\}.\]
\end{lem}
\begin{proof}
	First, we rewrite $\phi_{pcan}$ using $g(G) = |E(G)| - |V(G)| + 1$. Then,
	\begin{equation} \label{eqn:rewrite pcan}
		\phi_{pcan}(v) = \frac{|E(G)| + 1}{2|E(G)|}\cdot \deg(v) - 1.
	\end{equation}
	Using Lemma \ref{lem: subgraph_degrees} and equation \eqref{eqn:rewrite pcan}, for biconnected $W \subset V(G)$ we have
	\begin{align*}
		\phi_{pcan}(W) = \sum_{v \in W} \phi_{pcan}(v) = -|W| + \frac{|E(G)| + 1}{2|E(G)|}\left(|E(W, W^c)| + 2|E(G[W])|\right)
	\end{align*}
	Substituting this back into the definition of $\mathfrak{n}(\phi_{pcan})_W$, direct calculation gives
	\begin{align*}
		\mathfrak{n}(\phi_{pcan})_W =\left \lceil{\phi_{pcan}(W) - \frac{|E(W,W^c)|}{2}}\right \rceil = \left \lceil \frac{1}{2 |E(G)|}\sum_{v \in W} \deg(v)\right \rceil + |E(G[W])| - |W|
	\end{align*}
	By the Handshake Lemma, $\sum_{v \in W} \deg(v) < 2|E(G)|$, so $\mathfrak{n}(\phi_{pcan})_W = g(G[W]).$
\end{proof}
Theorem \ref{thm: bijection:vset-gbd} and Lemma \ref{lem:n_w(pcan)} yield the following corollary.
\begin{corollary} \label{cor:break-classical}
	For a graph $G$, the break divisors of a graph $G$, denoted $\mathcal{BD}_0(G)$ form a classical stability condition. In other words, $\mathcal{BD}_0(G) = \mathcal{BD}_{I_{\mathfrak{n}(\phi)}}$ for a generic numerical polarization $\phi$.
\end{corollary}
For the final result, we use a particular triangulating internal atlas. First, fix $q \in V(G)$ and consider each spanning tree $T$ of $G$ as a rooted tree with the root at $q$. For the internal atlas $\iA_q$, orient the internal edges of $T$ away from $q$. As shown in \cite{yuen2018} (see also \cite{ding2023framework}), this atlas is triangulating and corresponds to acyclic cocycle signature $\sigma_q^\ast$. Additionally, $(I_{G}^{\sigma^\ast_q} + (q))(T) = \vec{0}$ for every spanning tree $T$.
By Theorem \ref{thm: acyclic-regular}, $\left((\chi^\ast)^{-1}\circ \alpha^\ast\right)(\sigma_q^\ast)$ is a regular triangulation of the Lawrence polytope $\mathcal{P}^\ast$. Therefore, any other acyclic cocycle signature of $\M(G)$ is reachable from $\sigma^\ast_q$. 
\begin{thm}
	If $\sigma^\ast$ is an acyclic cocycle signature of the graphic matroid $\M(G)$, then $\mathcal{BD}_{I_G}^{\sigma^\ast}$ is a classical stability condition.
\end{thm}
\begin{proof}
	For some $q \in V(G)$, let $\sigma_q^\ast$ be the acyclic cocycle signature such that $I_G^{\sigma_q^\ast} = \vec{0}$. Let $\mathfrak{n} = \mathfrak{n}(\phi_{pcan})$. We induct on the number of cocycle flips $k$ needed to reach $\sigma^\ast$ from $\sigma_q^\ast$. If $\sigma_q^\ast = \sigma^\ast$, Corollary \ref{cor:break-classical} states that the break divisors are a classical stability condition, so $I_G^{\sigma_q^\ast} = I_{\mathfrak{n}}$. If $k=1$, let $C_1^\ast$ be the cocycle such that $\vec{C_1^\ast} \in \sigma_q^\ast$ and $-\vec{C_1^\ast} \in \sigma^\ast$. Let $W_1$ be the biconnected subset of $V(G)$ such that $C_1^\ast = E(W_1, W_1^c)$. Let $\phi = \phi_{pcan}$ and recall the wall and chamber decomposition of $V^d(G)$ induced by the hyperplane arrangement $\mathcal{A}_G^d$.

	Then, $\phi$ lives in the chamber cut out by the hyperplanes $\left\{\phi_W - \frac{|E(W,W^c)|}{2} = g(G[W]) \right\}$ for each biconnected $W \subset V(G)$. 

	Assume that $\vec{C_1^\ast} \in \sigma_q^\ast$ points towards $W_1$. Consider a generic polarization $\phi'$ in the chamber cut out by the same hyperplanes
	\[\left\{\phi_W - \frac{|E(W,W^c)|}{2} = g(G[W]) \right\}_{W \subset V(G), W \neq W_1, W_1^c} \]
	with the choice of new hyperplanes
	\[\left\{\phi_{W_1} - \frac{|E(W_1, W_1^c)|}{2} = g(G[W_1]) - 1 \right\} \bigcup \left \{\phi_{W_1^c} - \frac{|E(W_1, W_1^c)|}{2} = g(G[W_1^c]) + 1 \right\}.\]

Therefore, $\mathfrak{n}(\phi')_{w_1} = \mathfrak{n}(\phi)_{W_1} - 1$, $\mathfrak{n}(\phi')_{W_1^c} = \mathfrak{n}(\phi)_{W_1^c} + 1$, and $\mathfrak{n}(\phi')_W = \mathfrak{n}(\phi)_W$ for all other biconnected $W \subset V(G)$. With $\mathfrak{n}' = \mathfrak{n}(\phi')$, we can construct $I_{\mathfrak{n}'}$ as described in Theorem \ref{thm: bijection:vset-gbd}, and compare with $I_G^{\sigma^\ast}$. First, consider $T$ such that $C_1^\ast$ is not a fundamental cut of $T$.
Then, $\vec{T^\ast}$ does not change, so $I^{\sigma^\ast}(T)= I_G^{\sigma_q^\ast}(T)= \vec{0}$. Because $C_1^\ast$ is not a fundamental cut of $T$, there is no edge $e \in T$ such that $T\backslash e = T[W_1] \sqcup T[W_1^c]$, so $I_{\mathfrak{n}'}(T) = I_{\mathfrak{n}}(T) = \vec{0}$.

Now, consider $T$ such that $C_1^\ast$ is a fundamental cut of $T$. Let $\vec{T^\ast} \in \iA_{\sigma_q^\ast}$ and $\vec{T^{\ast'}} \in \iA_{\sigma^\ast}$, and let $e$ be the edge such that $C^\ast(T, e) = C_1^\ast$. Therefore, $\vec{e} \in \vec{T^\ast}$ and $-\vec{e} \in \vec{T^{\ast'}}$. If $\vec{e} = (v,v')$, we then have $I_G^{\sigma^\ast}(T)(v) = I_G^{\sigma_q^\ast}(T)(v) - 1$ and $I_G^{\sigma^\ast}(T)(v') = I_G^{\sigma_q^\ast}(T)(v') + 1$, since the indegree of $v$ decreases by $1$ and the indegree of $v'$ increases by $1$. Equivalently by construction of $\mathfrak{n}' = \mathfrak{n}(\phi')$, $	I_{\mathfrak{n'}}(T)_{W_1}=  I_{\mathfrak{n}}(T)_{W_1} - 1$ and $		I_{\mathfrak{n'}}(T)_{W_1^c} = I_{\mathfrak{n}}(T)_{W_1^c} + 1$.

The uniqueness of $I_{\mathfrak{n}'}$ guarantees $I_{\mathfrak{n}'} = I_G^{\sigma^\ast}$, as desired. Therefore, $\mathcal{BD}_{I_G^{\sigma^\ast}}(G)$ is a classical stability condition.

Now, assume the inductive hypothesis holds for all $0 \leq k \leq K$, and assume $\sigma^\ast$ is reachable from $\sigma_q^\ast$ through a sequence of $K+1$ cocycle flips. Let $C_1^\ast,...,C_{K+1}^\ast$ be the sequence of cocycles flipped, and assume $\sigma_K^\ast$ is the acyclic cocycle signature reached by the sequence $C_1^\ast,...,C_K^\ast$ of cocycle flips. By the inductive hypothesis, $\mathcal{BD}_{I_G^{\sigma_K^\ast}}(G)$ is a classical stability condition. Assume that $W_{K+1} \subset V(G)$ such that $C_{K+1}^\ast = E(W_{K+1}, W_{K+1}^c)$ and $\vec{C_{K+1}^\ast}$ is oriented towards $W_{K+1}$. If $\phi$ is the generic numerical polarization such that $I_{\mathfrak{n}(\phi)} = I_G^{\sigma_K^\ast}$, we can then select a generic numerical polarization $\phi'$ such that $\mathfrak{n}(\phi')_{W_{K+1}} = \mathfrak{n}(\phi)_{W_{K+1}} - 1$, $\mathfrak{n}(\phi')_{W_{K+1}^c} = \mathfrak{n}(\phi)_{W_{K+1}^c} + 1$, and $\mathfrak{n}(\phi')_{W} = \mathfrak{n}(\phi)_{W}$ for all other biconnected $W \subset V(G)$. By following the same procedure as the $k=1$ case, we see that $I_{\mathfrak{n}(\phi')} = I_G^{\sigma^\ast}$, as desired.  
\end{proof}

\end{document}